# Linear models of dynamic optimization with linear constraints

Somdeb Lahiri

ORCID: https://orcid.org/0000-0002-5247-3497

(Formerly with) PD Energy University, Gandhinagar (EU-G), India.

March 7, 2025.

This version: April 1, 2025.

## Abstract

We introduce a model of infinite horizon linear dynamic optimization with linear constraints and obtain results concerning feasibility of trajectories and optimal solutions necessarily satisfying conditions that resemble the "Euler condition" and "transversality condition". We prove results about optimal trajectories of strictly alternating problems, eventually conclusive problems, strongly eventually conclusive problems and two-phase problems.

**Keywords:** linear, dynamic optimization, optimal value function, dynamic programming, optimal decision rule, Euler condition, transversality condition, strictly alternating problems, eventually conclusive problems, strongly eventually conclusive problems, two-phase problems.

**AMS Subject Classification:** 90C05, 90C39, 90C46

**JEL Codes:** C44, C61.

## 1. Introduction:

One of the earliest works on infinite linear programming, i.e., linear programming with countably infinite variables and countably infinite linear constraints, is Evers (1973). Motivated by the models considered in Mitra (2000), Sorger (2015) and Lahiri (2025), we consider in the discussion that follows, a special case of the model presented in Evers (1973).

As in Mitra (2000), let $X = [0, b] \subset \mathbb{R}$, with $b > 0$ denote the **set of available alternatives**. With $\mathbb{N}$ denoting the set of natural number (i.e., the set of strictly positive integers) let $\mathbb{N}^0$ denote $\mathbb{N} \cup \{0\}$, i.e., the set of non-negative integers. Time is measured in discrete periods $t \in \mathbb{N}^0$. At each time 't' an alternative is chosen, and the chosen alternative is denoted by $x_t \in X$.

A very general form of the typical dynamic optimization problem along the lines defined in section 5.1 of Sorger (2015) is the following:

Given $x \in X$: Maximize $\sum_{t=0}^{\infty} u_t(x_t, x_{t+1})$, subject to the infinite sequence $<x_t | t \in \mathbb{N}^0>$ satisfying the constraints, $(x_t, x_{t+1}) \in \Omega_t$ for all $t \in \mathbb{N}^0$ and $x_0 = x \in X$, where for all $t \in \mathbb{N}^0$, $u_t$: $X \times X \to \mathbb{R}$ is the **utility function at time-period** t, $\Omega_t \subset X \times X$ is the **two-period constraint set at time-period** t, and $x \in X$ is the **initial choice**.

There is a voluminous literature devoted to the study of such models, particularly under the assumption that there exists a $\delta \in (0, 1)$ and $u: X \times X \to \mathbb{R}$ such that for all $t \in \mathbb{N}^0$ and $(x, y) \in X \times X$, $u_t(x, y) = \delta^t u(x, y)$ and $\Omega_t = \Omega_0$ for all $t \in \mathbb{N}^0$. Under such circumstances, time dependence is entirely confined to a discount factor that applies to a stationary utility function. This literature has been exhaustively surveyed in both Mitra (2000) and Sorger (2015).

In Lahiri (2025), we discuss a framework for the study of dynamic optimization problems in which the objective function is linear allowing the two-period constraint set to vary with time. The coefficients of the linear objective function are allowed to vary with time. The major theoretical results in Lahiri (2025) are the corner stone of the considerably general, and yet, a special case of the former, that we discuss here. In this paper, we study the special case where each time-dependent two-period constraint set is defined by a linear inequality, whose parameters are allowed to vary with time. One could call this a model of "linear programming in infinite dimensional spaces", but that – as pointed out in Lahiri (2025)- would be a misnomer for what is essentially a type of "asymptotic analysis".

With the two period constraint sets being defined by linear inequalities, in section 3, we prove simple properties that "feasible trajectories" of linear dynamic optimization problems with linear constraints satisfy. The names that we use for some of the properties, may not be the most appropriate ones, and we invite suggestions for improvement of terminology. In a subsequent section, we state results about optimal solutions that are inherited from the more general framework discussed in Lahiri (2025). In section 5, we use the theory of linear programming to obtain necessary conditions for optimality, that "resemble" the Euler condition and transversality condition. It has been shown in Lahiri (2025), that these two conditions, in their exact form, are sufficient for optimality. However, the necessary conditions we obtain here, are not the same, though they are very similar to the sufficient conditions. We also obtain necessary conditions for optimality, that rely on the first proposition about trajectories that satisfy the linear constraints, the latter being the first proposition of the paper.

Subcategories of problems we discuss are strictly alternating problems, eventually conclusive problems and two-phase problems. Strictly alternating problems are those whose coefficients in the objective function, beginning with the first time-period are non-zero and alternate in sign. Eventually conclusive problems are those whose coefficients in the objective function from a certain time-period onwards, are all non-positive. Two-phase problems are a further subcategory of eventually conclusive problems. For strictly alternating problems, with a mild additional requirement, we show that the optimal trajectory keeps oscillating between the upper bound and the lower bound. For eventually conclusive problems, there exists an optimal trajectory that always assumes the value zero after a finite number of time periods. Further, for problems that satisfy a stronger version of the eventually conclusive property, all optimal trajectories are such that they always assume the value zero after a finite number of time periods. For two phase problems, under a non-negativity assumption for both parameters defining the two-period constraint set during the first phase, in the initial phase the optimal trajectory adheres to the upper bound and eventually remains constant at zero.

Our verbal discussion of the results is far from adequate and those with an appetite for the exact results are cordially invited to proceed with their investigation of what follows.

**2. The Framework of Analysis:**

Let $\mathbb{R}$ denote the set of real numbers, $\mathbb{R}_+$ the set of non-negative real numbers and $\mathbb{R}_{++}$ the set of strictly positive real numbers. For any non-empty set A and $n \in \mathbb{N}$, we will represent the set of all n-tuples with coordinates in A by $A^n$. $A^n$ is the n-fold Cartesian product of A, i.e., the set of all functions from $\{1, \ldots, n\}$ to A.

As in section 1, let $X = [0, b] \subset \mathbb{R}$, with $b \in \mathbb{R}_{++}$ denote the **set of available alternatives**. Time is measured in discrete periods $t \in \mathbb{N}^0$. At each time 't' an alternative is chosen, and the chosen alternative is denoted by $x_t \in X$.

Let $<c_t | t \in \mathbb{N}^0>$ be a sequence in $\mathbb{R}_+$ and let $<a_t | t \in \mathbb{N}^0>$ be a sequence in $\mathbb{R}$ such that $\sum_{t=0}^{\infty} |c_t| < +\infty$, $\sum_{t=0}^{\infty} |a_t| < +\infty$ and for all $t \in \mathbb{N}^0$, $c_t + a_t b \geq 0$.

For all $x \in X = [0, b]$, $x = \frac{x}{b}b + (1-\frac{x}{b})0 = \frac{x}{b}b$.

Thus, for all $(x, t) \in X \times \mathbb{N}^0$: $c_t + a_t x = c_t + a_t \frac{x}{b}b = \frac{x}{b}(c_t + a_t b) + (1-\frac{x}{b})c_t \geq 0$, since $\frac{x}{b} \in [0, 1]$ and by hypothesis, both $c_t$ and $c_t + a_t b$ are non-negative.

For $t \in \mathbb{N}^0$, let $\Omega_t = \{(x, y) \in X \times X | y \leq c_t + a_t x\} = \{(x, y) | y \leq c_t + a_t x, x \in [0, b], y \in [0, b]\}$.

For $t \in \mathbb{N}^0$, $\Omega_t$ is the **two-period linearly constrained set at time-period** t.

It is precisely, this assumption of linearity in the definition of the "two-period constraint sets", that makes the model discussed here, a special case of the model discussed in Lahiri (2025).

For instance if for some $t \in \mathbb{N}^0$, $c_t = 0$, then for the same t it must be the case that $\Omega_t(0) = \{0\}$.

Note that for all $t \in \mathbb{N}^0$, $\Omega_t$ is a non-empty, closed and bounded subset of $X \times X$ and for each $x \in X$, the set $\Omega_t(x)$ **which is defined as** $\{y \in [0, b]| y \leq c_t + a_t x\} = [0, \min\{c_t + a_t x, b\}]$ is a non-empty, closed and bounded interval in X, though the interval $\Omega_t(x)$ may be a singleton (i.e., degenerate).

For $(x, t) \in X \times \mathbb{N}^0$, the set $\Omega_t(x)$ is said to be **the transition set from** x **at time-period** t.

For $x \in X$, let $\mathcal{F}(x) = \{<x_t | t \in \mathbb{N}^0> | x_{t+1} \in \Omega_t(x_t), t \in \mathbb{N}^0, x_0 = x\}$.

We will (whenever necessary) refer to an infinite sequence $<x_t | t \in \mathbb{N}^0> \in \mathcal{F}(x)$ as a **trajectory starting at (from)** x.

For $x \in X$ and $T \in \mathbb{N}^0$, let $\mathcal{F}^T(x) = \{<x_t | t \geq T> | x_{t+1} \in \Omega_t(x_t) \text{ for all } t \geq T, x_T = x\}$.

Clearly, $\mathcal{F}^T(x)$ is a convex set for all $(x, T) \in X \times \mathbb{N}^0$. Further, for all $x \in X$, $\mathcal{F}^0(x) = \mathcal{F}(x)$.

If $<x_t | t \geq T> \in \mathcal{F}^T(x)$ then we say that $<x_t | t \geq T>$ is a **trajectory starting at (from)** x **at time-period** T.

**Note 1:** For all $(x, T) \in X \times \mathbb{N}^0$, let $<x_t | t \geq T> \in \mathcal{F}^T(x)$ and $<y_t | t \geq T>$ be such that for some $T^0 \geq T$, $y_t = x_t$ for all $t \in \{T, \ldots, T^0\}$, $y_t = 0$ for all $t > T^0$. Clearly, $<y_t | t \geq T> \in \mathcal{F}^T(x)$.

Let $<p^{(t)} | t \in \mathbb{N}^0>$ be a sequence in $\mathbb{R}$ satisfying $\sum_{t=0}^{\infty} |p^{(t)}| < +\infty$

**Note 2:** $\sum_{t=0}^{\infty} |p^{(t)}| < +\infty$ **implies** $\lim_{t \to \infty} |p^{(t)}| = 0$ and hence for all sequence $<x_t| t \in \mathbb{N}^0>$ with $x_t \in X$ for all $t \in \mathbb{N}^0$, it must be the case that $\lim_{t \to \infty} |p^{(t)} x_t| = 0$.

We shall be concerned with is the following problem:

Given $x \in X$: Maximize $\sum_{t=0}^{\infty} p^{(t)} x_t$, subject to $x_{t+1} \in \Omega_t(x_t) = [0, \min\{c_t + a_t x, b\}]$ for all $t \in \mathbb{N}^0$, $x_0 = x$.

We shall refer to this problem as the **linear dynamic optimization with linear constraints (LDO-LC) problem** and represent it as $<(p^{(t)}, c_t, a_t)| t \in \mathbb{N}^0>$.

**Example 1: The linear cake eating problem:** In this case $a_t = 1$ and $c_t = 0$ for all $t \in \mathbb{N}^0$.

**Example 2: The linear optimal wealth accumulation problem:** $a_t \geq 0$ for all $t \in \mathbb{N}^0$ and there exists a sequence $<u_t| t \in \mathbb{N}^0> \mathbb{R}_+$ such that $\sum_{t=0}^{+\infty} u_t c_t \in \mathbb{R}$, $p_0 = u_0 a_0$, $p_t = u_t a_t - u_{t-1}$ for all $t \in \mathbb{N}$.

> For what follows we assume that $<((p^{(t)}, c_t, a_t)| t \in \mathbb{N}^0>$ is a given LDO-LC problem. As and when necessary, we will impose additional assumptions on this LDO problem.

## 3. Some Preliminary results about trajectories:

For any $(x, t) \in X \times \mathbb{N}^0$, it is easy to see that the following holds:

(i) $a_t > 0$ implies $c_t \leq c_t + a_t x \leq c_t + a_t b$, with at least one of the two-inequalities being strict.

(ii) $a_t = 0$ implies $c_t = c_t + a_t x = c_t + a_t b$.

(iii) $a_t < 0$ implies $c_t \geq c_t + a_t x \geq c_t + a_t b$, with at least one of the two-inequalities being strict

The following is a period-by-period characterization of trajectories that satisfy the constraints of the LDO-LC problem.

**Proposition 1:** Let $(x, T) \in X \times \mathbb{N}^0$ and let $<x_t| t \geq T> \in \mathcal{F}^T(x)$. Then for all $t > T$:

(i) $a_t > 0$ implies $x_t \in [\max\{0, \frac{x_{t+1} - c_t}{a_t}\}, \min\{c_{t-1} + a_{t-1} x_{t-1}, b\}]$.

(ii) $a_t = 0$ implies $x_t \in [0, \min\{c_{t-1} + a_{t-1} x_{t-1}, b\}]$.

(iii) $a_t < 0$ implies $x_t \in [0, \min\{\frac{x_{t+1} - c_t}{a_t}, \min\{c_{t-1} + a_{t-1} x_{t-1}, b\}\}]$.

**Proof:** Follows from noting that $<x_t| t \geq T> \in \mathcal{F}^T(x)$ implies that for all $t > T$, $x_t$ must satisfy $x \leq \min\{c_{t-1} + a_{t-1} x_{t-1}, b\}$, $a_t x \geq -c_t + x_{t+1}$, and $x \geq 0$. Q. E. D.

Suppose the sequence $<a_t| t \in \mathbb{N}^0>$ satisfies $a_t \geq 0$ for all $t \in \mathbb{N}^0$. Then, for any two sequences $<z_t| t \geq T>$ and $<w_t| t \geq T>$ satisfying $z_{t+1} = c_t + a_t z_t$ and $w_{t+1} = c_t + a_t w_t$, it is the case that [$z_t - w_t \geq 0$ if and only if $z_{t+1} - w_{t+1} \geq 0$]. On the other hand, if $a_t \leq 0$ for all $t \in \mathbb{N}^0$, then, for any two sequences $<z_t| t \geq T>$ and $<w_t| t \geq T>$ satisfying $z_{t+1} = c_t + a_t z_t$ and $z_{t+1} = c_t + a_t z_t$, it is the case that [$z_t - w_t \leq 0$ if and only if $z_{t+1} - w_{t+1} \geq 0$].

The following proposition follows from the observation above.

**Proposition 2:** Let $\langle (p^{(t)}, c_t, a_t) | t \in \mathbb{N}^0 \rangle$ be an LDO-LC problem with $a_t \leq 0$ for all $t \in \mathbb{N}^0$. Then the problem satisfies the following "**free disposability property**": For any $T \in \mathbb{N}^0$, and $y \in X$, $[\langle x_t | t \geq T \rangle \in \mathcal{F}^T(y)$, and $\langle y_t | t \geq T \rangle$ satisfies $y_T = x_T = y$, $x_t \geq y_t \geq 0$ for all $t > T]$ implies $[\langle y_t | t \geq T \rangle \in \mathcal{F}^T(y)]$.

**Proof:** For $t = T+1$, $0 \leq y_{T+1} \leq x_{T+1} \in [0, \min\{c_T + a_T y, b\}]$ and thus $y_{T+1} \in \Omega_T(y)$. For $t > T+1$ suppose $x_\tau \geq y_\tau \geq 0$ and $y_\tau \in \Omega_{\tau-1}(y_{\tau-1})$ for all $\tau = T+1, \ldots, t-1$.
Clearly, $x_t \in [0, \min\{c_{t-1} + a_{t-1} x_{t-1}, b\}]$ and $0 \leq y_t \leq x_t \leq \min\{c_{t-1} + a_{t-1} x_{t-1}, b\} \leq \min\{c_{t-1} + a_{t-1} y_{t-1}, b\}$, since $y_{t-1} \leq x_{t-1}$ and $a_{t-1} \leq 0$.

Thus, $y_t \in \Omega_t(y_{t-1})$.

By a standard induction argument, it follows that $\langle y_t | t \geq T \rangle \in \mathcal{F}^T(y)$. Q. E. D.

A proposition concerning different initial points but the same trajectory thereafter is the following.

**Proposition 3:** Let $\langle (p^{(t)}, c_t, a_t) | t \in \mathbb{N}^0 \rangle$ be an LDO-LC problem and suppose that for some $T \in \mathbb{N}^0$, it is the case that $a_T \geq 0$. Given $y > x$, let $\langle x_t | t \geq T \rangle \in \mathcal{F}^T(x)$, and let $\langle y_t | t \geq T \rangle$ be such that $y_T = y$, $y_t = x_t$ for all $t > T$. Then, $\langle y_t | t \geq T \rangle \in \mathcal{F}^T(y)$.

**Proof:** Since $y > x$ and $a_T \geq 0$, $\min\{c_T + a_T x, b\} \leq \min\{c_T + a_T y, b\}$. Thus, $x_{T+1} \in [0, \min\{c_T + a_T x, b\}] \subset [0, \min\{c_T + a_T y, b\}]$. Since $y_{T+1} = x_{T+1}$, it follows that $y_{T+1} \in [0, \min\{c_T + a_T y, b\}]$.

For $t > T+1$, $y_\tau = x_\tau$ for all $\tau = T+1, \ldots, t-1$, $x_t \in [0, \min\{c_{t-1} + a_{t-1} x_{t-1}, b\}]$ implies $y_t \in [0, \min\{c_{t-1} + a_{t-1} y_{t-1}, b\}]$.

Thus, $\langle y_t | t \geq T \rangle \in \mathcal{F}^T(y)$. Q.E.D.

**4. Results about optimal trajectories inherited from the general framework:**

For $(x, T) \in X \times \mathbb{N}^0$, consider the following problem: Maximize $\sum_{t=T}^{\infty} p^{(t)} x_t$, subject to $x_{t+1} \in \Omega_t(x_t) = [0, \min\{c_t + a_t x, b\}]$ for all $t \geq T$, $x_T = x$.

For $(x, T) \in X \times \mathbb{N}^0$, let $\mathcal{S}^T(x) = \underset{\langle x_t | t \geq T \rangle \in \mathcal{F}^T(x)}{\operatorname{argmax}} \sum_{t=T}^{\infty} p^{(t)} x_t$.

For $x \in X$, we will write $\mathcal{S}(x)$ for $\mathcal{S}^0(x)$.

Once again, it is quite straightforward to verify, that for all $(x, T) \in X \times \mathbb{N}^0$, $\mathcal{S}^T(x)$ is a convex set.

For $(x, T) \in X \times \mathbb{N}^0$, we refer to $\langle x_t | t \geq T \rangle \in \mathcal{S}^T(x)$ if $T > 0$ as an **optimal trajectory beginning at (from) x at time period** T, and if $T = 0$, then simply as an **optimal trajectory beginning at (from) x.**

By proposition 4.1 in Lahiri (2025) we have the following result.

**Proposition 4:** For $(x, T) \in X \times \mathbb{N}^0$, $\mathcal{S}^T(x) \neq \phi$.

By proposition 3, the following functions are well defined.

For $T \in \mathbb{N}^0$, $V^T: X \to \mathbb{R}$ is defined thus: for all $x \in X$, $V^T(x) = \sum_{t=T}^{\infty} p^{(t)} x_t$ for $<x_t | t \geq T> \in \mathcal{S}^T(x)$.

We will denote $V^0$ by $V$.

For $T \in \mathbb{N}$, $V^T$ is referred to as the **optimal value function for period** $T$, and $V$ is referred to as the **optimal value function**.

By propositions 5.1 and 5.2 in Lahiri (2025) we have the following result.

**Proposition 5:** (i) The optimal value function V, is concave and continuous on X.

(ii) For all $T \in \mathbb{N}^0$, $V^T$ is concave and continuous on X and satisfies the following ***functional equation of dynamic programming***: For all $T \in \mathbb{N}^0$, $z \in X$ and $<x_t | t \geq T> \in \mathcal{S}^T(z)$: $V^T(z) = p^{(T)}z + V^{T+1}(x_{T+1}) = p^{(T)}z + \max_{y \in \Omega_T(z)} \{V^{T+1}(y)\}$.

(iii) For all $x \in X$: $<x_t | t \in \mathbb{N}^0> \in \mathcal{S}(x)$ if and only if $<x_t | t \in \mathbb{N}^0> \in \mathcal{F}(x)$ and for all $T \in \mathbb{N}^0$ it is the case that $V^T(x_T) = p^{(T)} x_T + V^{T+1}(x_{T+1})$.

## 5. Results about optimal trajectories when constraints are linear:

A result about necessary conditions for a trajectory to be an optimal trajectory whose proof uses the duality theorem of linear programming is the following.

**Proposition 6:** For $(x, T) \in X \times \mathbb{N}^0$ let $<y_t | t \geq T> \in \mathcal{S}^T(x)$. Then there exists an infinite sequence $<(\lambda_t, \mu_t, \gamma_t)| t > T>$ in $\mathbb{R}^3_+$ such that for all $t > T$ and $y \in \Omega_{t-1}(y_{t-1})$ satisfying $(y, y_{t+1}) \in \Omega_t$: (i) $\lambda_t y \leq \lambda_t c_{t-1} + \lambda_t a_{t-1} y_{t-1} = \lambda_t y_t$; (ii) $-\mu_t a_t y \leq \mu_t c_t - \mu_t y_{t+1} = -\mu_t a_t y_t$; (iii) $\gamma_t y \leq \gamma_t b = \gamma_t y_t$; (iv) $\lambda_t - a_t \mu_t + \gamma_t \geq p^{(t)}$; and (v) $\lambda_t y_t - a_t \mu_t y_t + \gamma_t y_t = p^{(t)} y_t$ Further, $\lim_{t \to \infty} (\lambda_t - \mu_t a_t + \gamma_t) y_t = 0$.

**Proof:** For $(x, T) \in X \times \mathbb{N}^0$ let $<y_t | t \geq T> \in \mathcal{S}^T(x)$. Then for $t > T$, $y_t$ solves the following linear programming problem:

Maximize $\sum_{t=T}^{t-1} p^{(\tau)} y_\tau + p^{(t)} y + \sum_{\tau=t+1}^{\infty} p^{(\tau)} y_\tau$, subject to $y \leq c_{t-1} + a_{t-1} y_{t-1}$, $-a_t y \leq c_t - y_{t+1}$, $y \leq b$, $y \geq 0$.

The above maximization problem is equivalent to the following denoted LP:

Maximize $p^{(t)} y$ subject to $y \leq c_{t-1} + a_{t-1} y_{t-1}$, $-a_t y \leq c_t - y_{t+1}$, $y \leq b$, $y \geq 0$.

$y_t$ solves this problem <u>if and only if</u> there exists $\lambda_t, \mu_t, \gamma_t \geq 0$, such that along with $y_t$ the following holds:

(1) $y_t \leq c_{t-1} + a_{t-1} y_{t-1}$, $-a_t y_t \leq c_t - y_{t+1}$, $y_t \leq b$, $y_t \geq 0$.

(2) $\lambda_t c_{t-1} + \lambda_t a_{t-1} y_{t-1} = \lambda_t y_t$, $\mu_t c_t - \mu_t y_{t+1} = -\mu_t a_t y_t$, $\gamma_t y_t = \gamma_t b$.

(3) $\lambda_t - a_t \mu_t + \gamma_t \geq p^{(t)}$.

(4) $\lambda_t y_t - a_t \mu_t y_t + \gamma_t y_t = p^{(t)} y_t$.

For $y \in \Omega_{t-1}(x_{t-1})$ satisfying $\Omega_t(y, x_{t+1})$, clearly y satisfies $y \leq c_{t-1} + a_{t-1} y_{t-1}$, $-a_t y \leq c_t - y_{t+1}$, $y \leq b$, $y \geq 0$.

Since $\lambda_t$, $\mu_t$, $\gamma_t$ are all non-negative, along with these inequalities, we get that for $y \in \Omega_{t-1}(x_{t-1})$ satisfying $\Omega_t(y, x_{t+1})$, $\lambda_t y \leq \lambda_t c_{t-1} + \lambda_t a_{t-1} y_{t-1} = \lambda_t y_t$, $-\mu_t a_t y \leq \mu_t c_t - \mu_t y_{t+1} = -\mu_t a_t y_t$, $\gamma_t y \leq \gamma_t b = \gamma_t y_t$.

Along with (2), (3) and (4) above, we get (i) to (v) in the statement of the proposition.

From (4) and $\lim_{t \to \infty} p^{(t)} y_t = 0$, we get $\lim_{t \to \infty} (\lambda_t - \mu_t a_t + \gamma_t) y_t = 0$. Q.E.D.

**Note 3:** It is easy to see from (i) and (iii) of the proposition 6 (above) that if there exists $y > y_t$ satisfying $y \in \Omega_{t-1}(x_{t-1})$ and $(y, x_{t+1}) \in \Omega_t$, then the non-negativity of $\lambda_t$ and $\gamma_t$ implies that $\lambda_t = 0 = \gamma_t$.

An immediate consequence of proposition 6 is the following corollary.

**Corollary of Proposition 6:** For $(x, T) \in X \times \mathbb{N}^0$ let $\langle y_t | t \geq T \rangle \in \mathcal{S}^T(x)$. Then there exists an infinite sequence $\langle (\xi_t, \mu_t) | t > T \rangle$ in $\mathbb{R}_+^2$ such that for all $t > T$ and $y \in \Omega_{t-1}(x_{t-1})$ satisfying $(y, x_{t+1}) \in \Omega_t$: (i) $\xi_t y \leq \xi_t y_t$; (ii) $-\mu_t a_t y \leq \mu_t c_t - \mu_t y_{t+1} = -\mu_t a_t y_t$; (iii) $\xi_t - a_t \mu_t \geq p^{(t)}$; and (iii) $(\xi_t - \mu_t a_t) y_t = p^{(t)} y_t$. Further, $\lim_{t \to \infty} (\xi_t - \mu_t a_t) y_t = 0$.

**Proof:** The proof follows from proposition 6, by letting $\xi_t = \lambda_t + \gamma_t$. Q.E.D.

**Note 4:** The necessary conditions for optimality in Proposition 6 and its corollary, "resemble" the "Euler" condition and a version of the transversality condition, which are shown to be a sufficient condition for optimality in the more general context, in proposition 4.3 in Lahiri (2025).

We now provide a necessary condition for an optimal trajectory the proof of which is based on proposition 1.

**Proposition 7:** For $(x, T) \in X \times \mathbb{N}^0$ let $\langle x_t | t \geq T \rangle \in \mathcal{S}^T(x)$. Then, for $t > T$ the following conditions must be satisfied.

(i) $a_t > 0$ implies $[x_t = \min\{c_{t-1} + a_{t-1} x_{t-1}, b\}$ if $p^{(t)} > 0]$ & $[x_t = \max\{0, \frac{x_{t+1} - c_t}{a_t}\}$ if $p^{(t)} < 0]$.

(ii) $a_t = 0$ implies $[x_t = \min\{c_{t-1} + a_{t-1} x_{t-1}, b\}$ if $p^{(t)} > 0]$ & $[x_t = 0$ if $p^{(t)} < 0]$.

(iii) $a_t < 0$ implies $[x_t = \min\{\frac{x_{t+1} - c_t}{a_t}, \min\{c_{t-1} + a_{t-1} x_{t-1}, b\}\}$ if $p^{(t)} > 0]$ & $[x_t = 0$ if $p^{(t)} < 0]$.

**Proof:** Let $\langle x_t | t \geq T \rangle \in \mathcal{S}^T(x)$. Then as in the proof of proposition 6, it follows that for $t > T$, $y_t$ solves the following linear programming problem:

Maximize $p^{(t)} y$ subject to $y \leq c_{t-1} + a_{t-1} y_{t-1}$, $-a_t y \leq c_t - y_{t+1}$, $y \leq b$, $y \geq 0$.

Case 1: $a_t > 0$.

By (i) of proposition 1 we know that $x_t \in [\max\{0, \frac{x_{t+1} - c_t}{a_t}\}, \min\{c_{t-1} + a_{t-1} x_{t-1}, b\}]$.

If $p^{(t)} > 0$, then $x_t$ solves LP if and only if $x_t = \min\{c_{t-1} + a_{t-1} x_{t-1}, b\}$.

If $p^{(t)} < 0$, then $x_t$ solves LP if and only if $x_t = \max\{0, \frac{x_{t+1} - c_t}{a_t}\}$.

Case 2: $a_t = 0$.

By (ii) of proposition 1 we know that $x_t \in [0, \min\{c_{t-1} + a_{t-1}x_{t-1}, b\}]$.

If $p^{(t)} > 0$, then $x_t$ solves LP <u>if and only if</u> $x_t = \min\{c_{t-1} + a_{t-1}x_{t-1}, b\}$.

If $p^{(t)} < 0$, then $x_t$ solves LP <u>if and only if</u> $x_t = 0$.

Case 3: $a_t < 0$.

By (iii) of proposition 1 we know that $x_t \in [0, \min\{\frac{x_{t+1}-c_t}{a_t}, \min\{c_{t-1} + a_{t-1}x_{t-1}, b\}\}]$.

If $p^{(t)} > 0$, then $x_t$ solves LP <u>if and only if</u> $x_t = \min\{\frac{x_{t+1}-c_t}{a_t}, \min\{c_{t-1} + a_{t-1}x_{t-1}, b\}\}$.

If $p^{(t)} < 0$, then $x_t$ solves LP <u>if and only if</u> $x_t = 0$. Q.E.D.

## 6. Some interesting optimal solutions:

An LDO-LC problem $<(p^{(t)}, c_t, a_t)| t \in \mathbb{N}^0>$ is said to be a **strictly alternating LDO-LC problem** if $p^{(1)} \neq 0$ and for all $t \in \mathbb{N}$: (i) [$p^{(t)} \geq 0$ implies $p^{(t+1)} < 0$], (ii) [$p^{(t)} \leq 0$ implies $p^{(t+1)} > 0$].

One of the assumptions of our framework is that $\min\{c_t, c_t + a_t b\} \geq 0$ for all $t \in \mathbb{N}^0$.

The following result specifies the optimal decision rule for "a sub-class of" strictly alternating LDO-LC problems. We will require the above weak inequality to be strict for what follows.

**Proposition 8:** If $<(p^{(t)}, c_t, a_t)| t \in \mathbb{N}^0>$ is a strictly alternating LDO-LC problem satisfying $\min\{c_t, c_t + a_t b\} > 0$ for all $t \in \mathbb{N}^0$ then for $<y_t|t \in \mathbb{N}^0> \in \mathcal{S}^0(x)$:

(i) [$p^{(1)} > 0$ and $a_{2t} \leq 0$ for all $t \in \mathbb{N}$] implies [$y_{2t} = 0$ for all $t \in \mathbb{N}$, $y_1 = \min\{c_0 + a_0 x, b\}$ and $y_{2t-1} = \min\{c_{2t}, b\}$ for all $t > 1$].

(ii) [$p^{(1)} < 0$ and $a_{2t-1} \leq 0$ for all $t \in \mathbb{N}$] implies [for all $t \in \mathbb{N}$, $y_{2t} = \min\{c_{2t-1}, b\}$ and $y_{2t-1} = 0$].

**Proof:** Since $b > 0$ and by hypothesis $\min\{c_t, c_t + a_t b\} > 0$ for all $t \in \mathbb{N}^0$, it must be the case that $\min\{c_t + a_t x, b\} > 0$ for all $(x, t) \in X \times \mathbb{N}^0$.

(i) If $p^{(1)} > 0$ then $p^{(2t)} < 0$ and $p^{(2t-1)} > 0$ for all $t \in \mathbb{N}$.

$p^{(2t)} < 0$ and $a_{2t} \leq 0$ for all $t \in \mathbb{N}$ along with (ii) and (iii) of proposition 7 implies that $y_{2t} = 0$ for all $t \in \mathbb{N}$.

$p^{(2t-1)} > 0$, $a_{2t-1} \geq 0$ along with (i) combined with (ii) of proposition 7 implies $y_{2t-1} = \min\{c_{2t-2} + a_{2t-2}y_{2t-2}, b\}$.

$p^{(2t-1)} > 0$, $a_{2t-1} < 0$ and (iii) of proposition 7 implies $y_{2t-1} = \min\{\frac{-c_{2t-1}+y_{2t}}{a_{2t-1}}, \min\{c_{2t-2} + a_{2t-2}y_{2t-2}, b\}\} = \min\{\frac{-c_{2t-1}}{a_{2t-1}}, \min\{c_{2t-2} + a_{2t-2}y_{2t-2}, b\}\}$, since $y_{2t} = 0$.

Since $y_{2t} = 0 < c_{2t-1} + a_{2t-1}y_{2t-1}$, it is <u>not possible</u> that $y_{2t-1} = \frac{-c_{2t-1}+y_{2t}}{a_{2t-1}} = \frac{-c_{2t-1}}{a_{2t-1}}$.

Thus, $y_{2t-1} = \min\{\frac{-c_{2t-1}}{a_{2t-1}}, \min\{c_{2t-2} + a_{2t-2}y_{2t-2}\}\} = \min\{c_{2t-2} + a_{2t-2}y_{2t-2}\}$.

Since $y_{2t-2} = 0$ for all $t \geq 2$, we get $y_1 = \min\{c_0 + a_0 x, b\}$ and $y_{2t-1} = \min\{c_{2t-2}, b\}$ for $t \geq 2$.

(ii) If $p^{(1)} < 0$ then $p^{(2t)} > 0$ and $p^{(2t-1)} < 0$ for all $t \in \mathbb{N}$.

Thus, $p^{(2t-1)} < 0$ and $a_{2t-1} \leq 0$ for all $t \in \mathbb{N}$ implies by (ii) and (iii) of proposition 7 that $y_{2t-1} = 0$ for all $t \in \mathbb{N}$.

Since $y_{2t-1} = 0$ and $p^{(2t)} > 0$, for all $t \in \mathbb{N}$, by (i) and (ii) of proposition 7, for all $t \in \mathbb{N}$, $y_{2t} = \min\{c_{2t-1}, b\}$, whenever $a_{2t} \geq 0$.

Whenever $a_{2t} < 0$, by (iii) of proposition 7, for all $t \in \mathbb{N}$, $y_{2t} = \min\{\frac{y_{2t+1}-c_{2t}}{a_{2t}}, \min\{c_{2t-1} + a_{2t-1}y_{2t-1}, b\}\} = \min\{\frac{-c_{2t}}{a_{2t}}, \min\{c_{2t-1}, b\}\},$.

Since $y_{2t+1} = 0$ and $c_{2t} + a_{2t}y_{2t} > 0$, it is <u>not possible</u> that $y_{2t} = \frac{y_{2t+1}-c_{2t}}{a_{2t}} = \frac{-c_{2t}}{a_{2t}}$.

Thus, it must be the case that for all $t \in \mathbb{N}$, $y_{2t} = \min\{c_{2t-1}, b\}$ whenever $a_{2t} < 0$.

Thus, for all $t \in \mathbb{N}$, $y_{2t} = \min\{c_{2t-1}, b\}$. Q.E.D.

An LDO-LC problem is said to be **eventually conclusive** if there exists $T \in \mathbb{N}$ such that $p^{(t)} \leq 0$ for all $t \geq T$.

An eventually conclusive LDO-LC is said to be **strongly eventually conclusive** if there exists $T \in \mathbb{N}$ such that $p^{(t)} < 0$ for all $t \geq T$.

**Proposition 9:** An eventually conclusive LDO-LC problem has at least one optimal solution $\langle z_t | t \in \mathbb{N}^0 \rangle \in \mathcal{S}^0(x)$ such that for some $T \in \mathbb{N}$, $z_t = 0$ for all $t \geq T$, and if $T > 1$, then $\langle z_t | t = 0, \ldots, T-1 \rangle$ solves: Maximize $\sum_{t=0}^{T-1} p^{(0)} x_t$ subject to $x_{t+1} \leq c_t + a_t x_t$, $x_t \geq 0$, $x_t \leq b$, for all $t = 0, \ldots, T-2$, $x_0 = x$.

If it strongly eventually conclusive then there exists $T \in \mathbb{N}$ such that for all $\langle z_t | t \in \mathbb{N}^0 \rangle \in \mathcal{S}^0(x)$ it is the case that $z_t = 0$ for all $t \geq T$ and if $T > 1$, then $\langle z_t | t = 0, \ldots, T-1 \rangle$ solves: Maximize $\sum_{t=0}^{T-1} p^{(0)} x_t$ subject to $x_{t+1} \leq c_t + a_t x_t$, $x_t \geq 0$, $x_t \leq b$, for all $t = 0, \ldots, T-2$, $x_0 = x$.

**Proof:** Suppose there exists $T \in \mathbb{N}$ such that $p^{(t)} \leq 0$ for all $t \geq T$. Let $\langle y_t | t \in \mathbb{N}^0 \rangle \in \mathcal{S}^0(x)$ and let $\langle z_t | t \in \mathbb{N}^0 \rangle$ be such that $z_t = y_t$ for $t < T$ and $z_t = 0$ for all $t \geq T$. Clearly, $z_T = 0 \leq y_T = c_{T-1} + a_{T-1}y_{T-1} = c_{T-1} + a_{T-1}z_{T-1}$ and for $t > T$, $z_t = 0 \leq c_{t-1} = c_{t-1} + a_{t-1}0 = c_{t-1} + a_{t-1}z_{t-1}$.

Thus, $\langle z_t | t \in \mathbb{N}^0 \rangle \in \mathcal{F}(x)$ and $\sum_{t=0}^{\infty} p^{(t)} z_t = \sum_{t=0}^{T-1} p^{(t)} y_t \geq \sum_{t=0}^{\infty} p^{(t)} y_t$, since $p^{(t)} \leq 0$ and $y_t \geq 0$ for all $t \geq T$.

Since, $\langle y_t | t \in \mathbb{N}^0 \rangle \in \mathcal{S}^0(x)$ it must be the case that $\langle z_t | t \in \mathbb{N}^0 \rangle \in \mathcal{S}^0(x)$.

Thus, if $T > 1$, then $\langle z_t | t = 0, 1, \ldots, T-1 \rangle$ solves the problem: Maximize $\sum_{t=0}^{T-1} p^{(0)} x_t$ subject to $x_{t+1} \leq c_t + a_t x_t$, $x_t \geq 0$, $x_t \leq b$, for all $t = 0, \ldots, T-2$, $x_0 = x$.

Further, if the LDO-LC problem is "strongly eventually conclusive" and $\{t \geq T | y_t > 0\} \neq \phi$, then $\sum_{t=0}^{\infty} p^{(t)} z_t = \sum_{t=0}^{T-1} p^{(t)} y_t > \sum_{t=0}^{\infty} p^{(t)} y_t$, contradicting $\langle y_t | t \in \mathbb{N}^0 \rangle \in \mathcal{S}^0(x)$.

Hence, it must be the case that $y_t = 0$ for all $t \geq T$.

Thus, $\langle y_t | t \in \mathbb{N}^0 \rangle \in \mathcal{S}^0(x)$ if and only if [$y_t = 0$ for all $t \geq T$ and if $T > 1$, then $\langle y_t | t = 0, 1, \ldots, T-1 \rangle$ solves the problem: Maximize $\sum_{t=0}^{T-1} p^{(0)} x_t$ subject to $x_{t+1} \leq c_t + a_t x_t$, $x_t \geq 0$, $x_t \leq b$, for all $t = 0, \ldots, T-2$, $x_0 = x$. Q.E.D.

An LDO-LC problem $\langle (p^{(t)}, c_t, a_t) | t \in \mathbb{N}^0 \rangle$ is said to be a **two-phase LDO-LC problem** if there exists $T^+ \in \mathbb{N}^0$, $T^- \in \mathbb{N}$ with $T^+ < T^-$ such that $T^+ = \max\{T | p^{(t)} > 0 \text{ for all } t \leq T\}$ and $T^- = \min\{T | p^{(t)} < 0 \text{ for all } t \geq T\}$.

Clearly a two-phase LDO-LC problem is strongly eventually conclusive.

For such problems we have the following proposition about optimal trajectories.

**Proposition 10:** Consider a two-phase LDO-LC problem satisfying $a_t \geq 0$ for all $t < T^+$ if $T^+ > 0$.

If $\langle y_t | t \in \mathbb{N}^0 \rangle \in \mathcal{S}^0(x)$, then $y_t = 0$ for all $t \geq T^-$ and if $T^+ \geq 1$ then for all $t = 0, \ldots, T^+-1$, $y_{t+1} = \min\{c_t + a_t y_t, b\}$.

**Proof:** Let $\langle y_t | t \in \mathbb{N}^0 \rangle \in \mathcal{S}^0(x)$. Since a two-phase LDO-LC is strongly eventually conclusive, by proposition 9 we know that $y_t = 0$ for all $t \geq T^-$.

Suppose $T^+ \geq 1$ and towards a contradiction suppose that for some $t \in \{0, \ldots, T^+-1\}$, $y_{t+1} < \min\{c_t + a_t y_t, b\}$.

Let $T = \min\{t | t \in \{0, \ldots, T^+-1\}, y_{t+1} < \min\{c_t + a_t y_t, b\}\}$. Let $\langle z_t | t \in \mathbb{N}^0 \rangle$ be such that $z_t = y_t$ for $t \neq T+1$, $z_{T+1} = \min\{c_T + a_T y_T, b\}$. Since $z_{T+1} > y_{T+1}$ and $a_{T+1} \geq 0$, $c_{T+1} + a_{T+1} z_{T+1} \geq c_{T+1} + a_{T+1} y_{T+1}$, so that $\min\{c_{T+1} + a_{T+1} z_{T+1}, b\} \geq \min\{c_{T+1} + a_{T+1} y_{T+1}, b\} \geq y_{T+2} = z_{T+2}$.

Thus, $\langle z_t | t \in \mathbb{N}^0 \rangle \in \mathcal{F}(x)$.

Further, $\sum_{t=0}^{\infty} p^{(t)} z_t = \sum_{t=0}^{\infty} p^{(t)} y_t + p^{(T+1)}(z_{T+1} - y_{T+1}) > \sum_{t=0}^{\infty} p^{(t)} y_t$, since $p^{(T+1)} > 0$ and $z_{T+1} - y_{T+1} > 0$.

This contradicts our assumption that $\langle y_t | t \in \mathbb{N}^0 \rangle \in \mathcal{S}^0(x)$.

Thus, if $T^+ \geq 1$ then for all $t = 0, \ldots, T^+-1$, $y_{t+1} = \min\{c_t + a_t y_t, b\}$. Q.E.D.